\documentclass[12pt]{amsart}
\usepackage{amsmath,amsthm,amsfonts,amscd,amssymb,eucal,latexsym}

\theoremstyle{definition}
\newtheorem{theorem}{Theorem}[section]

\newtheorem{lemma}[theorem]{Lemma}
\newtheorem{proposition}[theorem]{Proposition}
\newtheorem{remark}[theorem]{Remark}
\newtheorem{notation}[theorem]{Notation}

\newtheorem{definition}[theorem]{Definition}

\newcommand{\rank}{{\rm rank}}
\newcommand{\CPA}{{\rm CPA}}
\newcommand{\rcp}{{\rm rcp}}

\newcommand{\spn}{{\rm span}}

\newcommand{\spec}{{\rm spec}}
\newcommand{\op}{{\rm op}}
\newcommand{\Aut}{{\rm Aut}}

\allowdisplaybreaks

\begin{document}

\title[noncommutative toral automorphisms]{Positive Voiculescu-Brown
entropy in noncommutative toral automorphisms}

\author{David Kerr}
\author{Hanfeng Li}
\address{Dipartimento di Matematica, Universit\`{a} di Roma
``La Sapienza,'' P.le Aldo Moro, 2, 00185 Rome, Italy}
\email{kerr@mat.uniroma1.it}
\address{Department of Mathematics, University of Toronto, Toronto,
Ontario M5S 3G3, Canada}
\email{hli@fields.toronto.edu}
\date{March 7, 2003}

\begin{abstract}
We show that the Voiculescu-Brown entropy of a noncommutative toral
automorphism arising from a matrix
$S\in GL(d,\mathbb{Z})$ is at least half the value of the topological
entropy of the corresponding classical toral automorphism. We
also obtain some information concerning the
positivity of local Voiculescu-Brown entropy with respect to single unitaries.
In particular we show that
if $S$ has no roots of unity as eigenvalues then the
local Voiculescu-Brown entropy with respect to every product of
canonical unitaries is positive, and also that in the presence of
completely positive CNT entropy the unital version of local
Voiculescu-Brown entropy with respect to every non-scalar unitary is
positive.
\end{abstract}

\maketitle

\section{Introduction}

Let $\Theta = (\theta_{jk} )_{1\leq j,k \leq d}$ be a real
skew-symmetric $d\times d$ matrix. The noncommutative $d$-torus
$A_\Theta$ is defined as the universal $C^*$-algebra generated by
unitaries $u_1 , \dots , u_d$ subject to the relations
$$ u_j u_k = e^{2\pi i \theta_{jk}} u_k u_j $$
for all $1\leq j,k \leq d$ (see \cite{Rie} for a reference).
For any matrix $S= (s_{jk} )_{1\leq
j,k \leq d}$ in $GL(d,\mathbb{Z})$ there is an isomorphism $\alpha
: A_{S^t \Theta S} \to A_\Theta$ determined by
$$ \alpha_\Theta (u_j ) = u^{s_{1j}}_1 u^{s_{2j}}_2 \cdots u^{s_{dj}}_d $$
for each $j=1, \dots ,d$. Thus when $S^t \Theta S \equiv
\Theta\,\,\, (\text{\rm mod}\,\, M_d (\mathbb{Z}))$ we obtain an
automorphism of $A_\Theta$, which we denote by $\alpha_{S,\Theta}$
and refer to as a {\it noncommutative toral automorphism}. Note
that whenever $\alpha_{S,\Theta}$ exists so does
$\alpha_{S,-\Theta}$. These noncommutative analogues of toral
automorphisms were initially introduced in \cite{Wat} and
\cite{Bre} for $d=2$, in which case for any given $S\in
SL(2,\mathbb{Z})$ the automorphism $\alpha_{S,\Theta}$ is defined
for all $\Theta$. An indication of their significance from a
noncommutative geometry perspective is the fact that, for $d=2$,
if $\theta_{12}$ is an irrational number satisfying a generic
Diophantine property, then every diffeomorphism of $A_\Theta$
equipped with the canonical differential structure is a
composition of an inner automorphism by a smooth unitary in the
connected component of the unit, a noncommutative toral
automorphism, and an automorphism arising from the canonical
action of $\mathbb{T}^2$ \cite{Ell}. In an arbitrary dimension
$d$, if $\Theta$ is rational (i.e., the entries of $\Theta$ are
all rational) then the $C^*$-algebra $A_\Theta$ is homogeneous
(see \cite{OPT}) and when $\alpha_{S,\Theta}$ exists we recover
the corresponding classical toral automorphism at the level of the
pure state space upon restricting $\alpha_{S,\Theta}$ to the
center of $A_\Theta$, so that from the noncommutative point of
view it is the case of nonrational $\Theta$ that is of primary
interest. Unlike classical toral automorphisms, which in the
hyperbolic case have served as prototypes for such important
phenomena as hyperbolicity and structural stability (see
\cite{DGS,HK}), noncommutative toral automorphisms as a class have
remained somewhat mysterious, although in many cases much has been
ascertained from a measure-theoretic viewpoint
\cite{BNS,Nar,NT,Nar2,Nes}.

In \cite[Sect.\ 5]{Voi} it was shown that the Voiculescu-Brown
entropy of $\alpha_{S,\Theta}$ is bounded above by the topological
entropy of the corresponding classical toral automorphism, i.e.,
by $\sum_{| \lambda_i | > 1} \log | \lambda_i |$, where $\lambda_1
, \dots , \lambda_d$ are the eigenvalues of $S$ counted with
multiplicity (see, e.g., \cite[Sect.\ 24]{DGS} for a calculation
of the topological entropy of toral automorphisms). Except in the
case that the eigenvalues of $S$ all lie on the unit circle, this
does not resolve the basic question of whether the entropy is
positive or zero, i.e., of whether the system is ``chaotic'' or
``deterministic.'' When $\Theta$ is rational the Voiculescu-Brown
entropy can be seen to be $\sum_{| \lambda_i | > 1}\log |
\lambda_i |$ by restricting $\alpha_{S,\Theta}$ to the centre of
$A_\Theta$ and applying Propositions 4.4 and 4.8 of \cite{Voi}. In
\cite{Nes} Neshveyev showed that if $S$ has no roots of unity
then the entropic $K$-property (and in particular
positive CNT entropy and hence also positive Voiculescu-Brown
entropy by \cite[Prop.\ 4.6]{Voi}) follows from a summability
condition with respect to a $2$-cocycle $\mathbb{Z}^d
\times\mathbb{Z}^d \to\mathbb{T}$ in terms of which $A_\Theta$ can
be described. In the case $d=2$, when we are dealing with a
rotation $C^*$-algebra $A_\theta$ ($=A_\Theta$ for $\Theta =
(\theta_{jk})_{j,k=1,2}$ with $\theta_{12} = \theta$), if $S$ has
an eigenvalue $\lambda$ with $|\lambda| > 1$ then the set of
$\theta$ for which the CNT entropy with respect to the canonical
tracial state is positive has zero Lebesgue measure \cite{NT} and
it contains $\mathbb{Z} + 2\mathbb{Z}\lambda^2$ as a consequence
of the above-mentioned result of Neshveyev (see \cite{Nar,Nes}).

The first aim of this article, which we carry out in Section~\ref{S-pos},
is to show that, in an arbitrary dimension
$d$, if the eigenvalues $\lambda_1 , \dots , \lambda_d$ do not all
lie on the unit circle, then the
Voiculescu-Brown entropy of $\alpha_{S,\Theta}$ is at least
$\frac12 \sum_{| \lambda_i | > 1} \log | \lambda_i |$.
(This was effectively claimed in \cite{NT} for $d=2$ but the tensor
product argument given there is not correct.) In Section~\ref{S-posloc}
we apply a result from \cite{EID} to obtain some information
concerning the positivity of the local Voiculescu-Brown entropy of
$\alpha_{S,\Theta}$
with respect to products of canonical unitaries. We prove in particular
that if $S$ has no roots of unity as eigenvalues then
the local Voiculescu-Brown entropy of $\alpha_{S,\Theta}$ with respect
to any product of canonical unitaries is positive.
Finally, in Section~\ref{S-K} we show in the general unital setting that
completely positive CNT entropy of the von Neumann algebraic dynamical
system arising from a faithful invariant state
implies positivity of the unital version of local Voiculescu entropy
with respect to every non-scalar
unitary, and also that in the unital case this latter condition is
equivalent to a
noncommutative extension of the topological-dynamical notion of
completely positive entropy which we call
``completely positive Voiculescu-Brown entropy.''
We then apply these results to the above-mentioned noncommutative
toral automorphisms treated by Neshveyev in \cite{Nes}.
\medskip

\noindent{\it Acknowledgements.} D. Kerr was supported by the
Natural Sciences and Engineering Research Council of Canada. This
work was carried out during his stays at the University of Tokyo
and the University of Rome ``La Sapienza'' over the academic years
2001--2002 and 2002--2003, respectively. He thanks Yasuyuki
Kawahigashi at the University of Tokyo and Claudia Pinzari at the
University of Rome ``La Sapienza'' for their generous hospitality.
H. Li was supported jointly by the Mathematics Department of the
University of Toronto and NSERC Grant 8864-02 of George A.
Elliott. We thank George Elliott for valuable discussions. H.
Li would also like to thank Sergey Neshveyev for helpful
discussions.

\section{Positive Voiculescu entropy}\label{S-pos}

We begin by recalling the definition of Voiculescu-Brown entropy
\cite{Br}, which is based on completely positive approximation
(see \cite{Pau} for a reference on completely positive maps). Let
$A$ be an exact $C^*$-algebra and $\alpha$ an automorphism of $A$.
Let $\pi : A \to\mathcal{B} (\mathcal{H})$ be a faithful
$^*$-representation. For a finite set $\Omega\subseteq A$ and
$\delta > 0$ we denote by $\CPA (\pi , \Omega , \delta )$ the
collection of triples $(\varphi , \psi , B)$ where $B$ is a
finite-dimensional $C^*$-algebra and $\varphi : A\to B$ and $\psi
: B\to\mathcal{B} (\mathcal{H})$ are contractive completely
positive maps such that $\| (\psi\circ\varphi )(a) - \pi (a) \| <
\delta$ for all $a\in\Omega$, and we define $\rcp (\Omega , \delta
)$ to be the infimum of $\rank\, B$ over all $(\varphi , \psi ,
B)\in\CPA (\pi , \Omega , \delta )$, with rank referring to the
dimension of a maximal Abelian $C^*$-subalgebra. As the notation
indicates, $\rcp (\Omega , \delta )$ is independent of the
faithful $^*$-representation $\pi$, as shown in the proof of
Proposition 1.3 in \cite{Br}. We then set
\begin{align*}
ht(\alpha , \Omega , \delta ) &= \limsup_{n\to\infty}\frac1n \log\rcp
(\Omega\cup\alpha\Omega\cup\cdots\cup\alpha^{n-1}\Omega , \delta ) , \\
ht(\alpha , \Omega ) &= \sup_{\delta > 0}ht(\alpha , \Omega , \delta ) , \\
ht(\alpha ) &= \sup_\Omega ht(\alpha , \Omega )
\end{align*}
where the last supremum is taken over all finite sets
$\Omega\subseteq A$. The quantity $ht(\alpha )$ is a
$C^*$-dynamical invariant which we call the {\it Voiculescu-Brown
entropy} of $\alpha$.

We will also have occasion to use the unital
version of $ht(\alpha , \Omega )$ in Section~\ref{S-K}, which we will denote
by $ht'(\alpha , \Omega )$. This is defined in the case of unital $A$ by
using unital completely positive maps instead of general contractive
completely positive maps as above. When $1\not\in\Omega$ we may have
$ht(\alpha , \Omega )\neq ht'(\alpha , \Omega )$, but these quantities
do agree when
$1\in\Omega$ as the proof of Proposition 1.4 in \cite{Br} shows, so
that $ht(\alpha )$ may be alternatively obtained in the
unital case by taking the supremum of $ht'(\alpha , \Omega )$ over
all finite sets $\Omega\subseteq A$.

\begin{notation}\label{N-op}
For a $C^*$-algebra $A$ we denote by $A^{\op}$ the opposite algebra
(i.e., the $C^*$-algebra obtained from $A$ by reversing the multiplication),
and for each $a\in A$ we denote by $\tilde{a}$ the corresponding element in
$A^{\op}$.
\end{notation}

We would like to thank George Elliott for suggesting the proof we
give of the following lemma, which simplifies our original proof.

\begin{lemma}\label{L-op}
Let $A$ and $B$ be $C^*$-algebras and $\varphi : A\to B$ an $n$-positive
map. Let $\varphi^{\op} : A^{\op}\to B^{\op}$ be the induced linear map
given by $\varphi^{\op} (\tilde{a}) = \widetilde{\varphi (a)}$ for all
 $a\in A$. Then $\varphi^{\op}$ is also $n$-positive. In particular,
if $\varphi$ is completely positive then so is $\varphi^{\op}$.
\end{lemma}

\begin{proof}
Suppose that $\varphi$ is $n$-positive, i.e., the map
$\varphi\otimes id:A\otimes M_n(\mathbb{C})\to B\otimes
M_n(\mathbb{C})$ is positive. Since for every $C^*$-algebra $D$
the map $D\rightarrow D^{\op}$ given by $a\mapsto \tilde{a}$ for
all $a\in D$ is an order isomorphism, it follows that the map
$(\varphi\otimes id)^{\op} : (A\otimes M_n(\mathbb{C}))^{\op}\to
(B\otimes M_n(\mathbb{C}))^{\op}$ is positive. But this is just
the map $\varphi^{\op}\otimes id^{\op} : A^{\op}\otimes
(M_n(\mathbb{C}))^{\op}\to B^{\op}\otimes
(M_n(\mathbb{C}))^{\op}$. Take an isomorphism
$\beta:(M_n(\mathbb{C}))^{op}\to M_n(\mathbb{C})$. Then
$(id_{B^{\op}}\otimes \beta)\circ (\varphi^{\op}\otimes
id^{\op})\circ (id_{A^{\op}}\otimes
\beta)^{-1}=\varphi^{\op}\otimes (\beta \circ id^{\op}\circ
\beta^{-1})=\varphi^{\op}\otimes id:A^{\op}\otimes
M_n(\mathbb{C})\rightarrow B^{\op}\otimes M_n(\mathbb{C})$ is
positive, i.e., $\varphi^{\op}$ is $n$-positive.
\end{proof}

\begin{proposition}\label{P-op}
Let $\alpha$ be an automorphism of an exact $C^*$-algebra $A$ and
let $\alpha^{\op}$ be the induced automorphism of $A^{\op}$. Then
$ht(\alpha ) = ht(\alpha^{\op})$.
\end{proposition}

\begin{proof}
Given a faithful $^*$-representation $\pi : A\to\mathcal{B} (\mathcal{H})$
consider the induced injective $^*$-homomorphism $\pi^{\op} : A^{\op} \to
\mathcal{B} (\mathcal{H} )^{\op}$. Let $\tilde{\mathcal{H}} = \{ \tilde{x} :
x\in\mathcal{H} \}$ be the Hilbert space
conjugate to $\mathcal{H}$ with scalar multiplication $\bar{\lambda}
\tilde{x} = \widetilde{\lambda x}$ and inner product
$\left<\tilde{x}, \tilde{y}\right> = \left<y,x\right>$. Then we have a
natural identification of $\mathcal{B} (\mathcal{H} )^{\op}$ with
$\mathcal{B} (\tilde{\mathcal{H}})$ under which we can consider $\pi^{\op}$
as a $^*$-representation. It follows then by Lemma~\ref{L-op} that for any
finite set $\Omega\subseteq A$ and $\delta > 0$ we have
$\rcp (\tilde{\Omega} , \delta ) = \rcp (\Omega , \delta )$ where
$\tilde{\Omega} = \{ \tilde{a} : a\in\Omega \}$,
and so we conclude that $ht(\alpha ) = ht(\alpha^{\op})$.
\end{proof}

\begin{lemma}\label{L-circle}
Let $d\geq 2$ and let $\Theta$ be a real skew-symmetric $d\times
d$ matrix. Let $\beta$ be the canonical action of $\mathbb{T}^d$
on the noncommutative torus $A_{\Theta}$, and let $\alpha$ be an
automorphism of $A_{\Theta}$ such that $\alpha
\beta(\mathbb{T}^n)\alpha^{-1}=\beta(\mathbb{T}^n)$ in
$\Aut(A_{\Theta})$. Then
$ht(\beta_x\alpha)=ht(\alpha\beta_x)=ht(\alpha)$ for all $x\in
\mathbb{T}^n$.
\end{lemma}

\begin{proof} It suffices to show $ht(\beta_x\alpha)\leq ht(\alpha)$.
Let $\pi : A_\Theta \to \mathcal{B}(\mathcal{H})$ be a faithful
$^*$-representation of $A_\Theta$. For each $p\in \mathbb{Z}^d$
there is a unitary $u_p$ in $A_\Theta$ such that $\spn \{ u_p :
p\in \mathbb{Z}^d \}$ is dense in $A_\Theta$ and $\beta_x (u_p) =
\left< p, x \right> u_p$, where $\left< \cdot,
\cdot\right>:\mathbb{Z}^d\times \mathbb{T}^d\to \mathbb{T}$ is the
canonical pairing. For every $\omega\subseteq \mathbb{Z}^d$ we set
$U_\omega = \{ u_p : p\in\omega \}$. Since $U_{\mathbb{Z}^d}$ is
total in $A_\Theta$, by Proposition 2.6 of \cite{Br} $ht(\alpha)$
and $ht(\beta_x\alpha)$ are equal to the supremum of $ht(\alpha,
U_\omega )$ and $ht(\beta_x\alpha, U_\omega)$ over all finite sets
$\omega\subseteq\mathbb{Z}^d$ respectively. Thus it suffices to show
that $ht(\alpha, U_\omega ) \leq ht(\beta_x \alpha, U_\omega )$ for
every finite set $\omega\subseteq\mathbb{Z}^d$, and this will follow
once we show that
\begin{gather*}
\rcp (U_\omega \cup (\beta_x \alpha )(U_\omega )\cup\cdots\cup
(\beta_x \alpha )^{m-1}(U_\omega ), \delta ) \hspace*{25mm} \\
\hspace*{35mm} \leq \rcp
(U_\omega \cup\alpha (U_\omega )\cup\cdots\cup\alpha^{m-1} (U_\omega ),
\delta )
\end{gather*}
for any given finite set $\omega\subseteq\mathbb{Z}^d$, $m\in\mathbb{N}$,
and $\delta>0$. Suppose then that $(\varphi, \psi, B)$
is a triple in $\CPA (\pi, U_\omega\cup\alpha (U_\omega )\cup\cdots\cup
\alpha^{m-1} (U_\omega ), \delta )$ such that
$\rank (B) = \rcp (U_\omega \cup\alpha (U_\omega )\cup\cdots\cup
\alpha^{m-1} (U_\omega ), \delta)$. For each $j\in \mathbb{Z}_{\geq 0}$
there exists some $x(j)\in\mathbb{T}^d$ such that
$(\beta_x \alpha )^j=\alpha^j \beta_{x(j)}$.
Then $(\beta_x \alpha )^j(u_p ) = \left< p,x(j) \right> \alpha^j (u_p )$ for
every $p\in\mathbb{Z}^d$, and so
$$ \| (\psi\circ\varphi )((\beta_x \alpha )^j(u_p )) -
\pi ((\beta_x \alpha )^j (u_p )) \| = \| (\psi\circ\varphi )
(\alpha^j (u_p )) - \pi (\alpha^j (u_p )) \| < \delta $$ for all
$j=0, \dots , m-1$ and $p\in\omega$. Thus the triple $(\varphi, \psi, B)$
is also contained in $\CPA (\pi, U_\omega\cup (\beta_x \alpha )(U_\omega )
\cup\cdots\cup (\beta_x \alpha )^{m-1} (U_\omega ), \delta )$,
finishing the proof.
\end{proof}

\begin{remark}\label{R-tori}
(1) Let $A$ be any exact $C^*$-algebra with a sequence of
finite dimensional subspaces $V_1\subseteq V_2\subseteq \cdots$
such that $\bigcup_{j\in\mathbb{N}} V_j$ is dense in $A$, and let $G$ be a
subgroup of $\Aut (A)$ preserving every $V_j$. If $\alpha\in\Aut (A)$
satisfies $\alpha G\alpha^{-1}=G$ then
$ht(\beta \alpha ) = ht(\alpha \beta ) = ht(\alpha )$ for every
$\beta\in G$. The proof of Lemma~\ref{L-circle} applies with minor
modifications.

(2) It is easy to show that an automorphism $\alpha$ of $A_\Theta$
satisfies the hypothesis of Lemma~\ref{L-circle} if and only if
it is of the form $\alpha_{S,\Theta}\beta_x$ for some noncommutative
toral automorphism $\alpha_{S,\Theta}$ and $x\in \mathbb{T}^d$.
\end{remark}

\begin{lemma}\label{L-tensor}
Let $\alpha_{S,\Theta}$ be any noncommutative toral automorphism. Then
$$ ht(\alpha_{S,\Theta}) + ht(\alpha_{S,-\Theta}) \geq
\sum_{| \lambda_i | > 1} \log | \lambda_i | $$
where $\lambda_1 , \dots , \lambda_d$ are the eigenvalues of $S$
counted with multiplicity.
\end{lemma}

\begin{proof}
Since the Voiculescu-Brown entropy of an automorphism of a separable
commutative $C^*$-algebra agrees with the topological entropy of the induced
homeomorphism on the pure state space by \cite[Prop.\ 4.8]{Voi}, in the case
$\Theta = 0$ we have $ht(\alpha_{S,0} ) =
\sum_{| \lambda_i | > 1} \log | \lambda_i |$.
Consider now the tensor product $A_\Theta \otimes A_{-\Theta}$. Denoting
by $u_1 , \dots , u_d$ and $v_1 , \dots , v_d$ the canonical
unitaries of $A_\Theta$ and $A_{-\Theta}$, respectively, we have that
the unitaries $u_j \otimes v_j$ for $j=1, \dots ,d$ form
canonical generators for a copy $C$ of $C(\mathbb{T}^d )$. This can been
seen from the fact that they operate as shifts in different coordinate
directions on the Hilbert subspace
$$ \overline{\spn} \{ \pi_+ (u^{k_1}_1 \cdots u^{k_d}_d )
\xi_+ \otimes\pi_- (v^{k_1}_1 \cdots v^{k_d}_d )\xi_- : (k_1 ,
\dots , k_d ) \in\mathbb{Z}^d \} $$ (identified with $\ell^2
(\mathbb{Z}^d )$) with respect to the tensor product of the
canonical tracial state GNS representations $\pi_{\pm}$ of
$A_{\pm\Theta}$ with canonical cyclic vectors $\xi_{\pm}$. We
furthermore see that this identification of the $\alpha_{S,\Theta}
\otimes \alpha_{S,-\Theta}$-invariant $C^*$-subalgebra $C$ with
$C(\mathbb{T}^d )$ establishes a conjugacy between
$\alpha_{S,\Theta} \otimes \alpha_{S,-\Theta} \big| {}_C$ and
$\alpha_{S,0}$. The monotonicity and tensor product subadditivity
of Voiculescu-Brown entropy then yields
$$ \sum_{| \lambda_i | > 1} \log | \lambda_i | =
ht(\alpha_{S,0} ) = ht(\alpha_{S,\Theta} \otimes \alpha_{S,-\Theta}
\big| {}_C ) \leq ht(\alpha_{S,\Theta}) + ht(\alpha_{S,-\Theta}) . $$
\end{proof}

\begin{theorem}\label{T-half}
Let $\alpha_{S,\Theta}$ be any noncommutative toral automorphism. Then
$$ ht(\alpha_{S,\Theta}) =ht(\alpha_{S,-\Theta})\geq
\frac12 \sum_{| \lambda_i | > 1} \log | \lambda_i | $$ where
$\lambda_1 , \dots , \lambda_d$ are the eigenvalues of $S$ counted
with multiplicity.
\end{theorem}

\begin{proof}
Denoting by $u_1 , \dots , u_d$ and $v_1 , \dots , v_d$ the canonical
unitaries of $A_\Theta$ and $A_{-\Theta}$, respectively, we have an
isomorphism $A_\Theta \to A^{\op}_{-\Theta}$ given by
$u_j \mapsto \tilde{v_j}$. Identify $A_\Theta$ and
$A^{\op}_{-\Theta}$ via this isomorphism. By computing
$\alpha^{\op}_{S,-\Theta} \in\Aut (A_\Theta )$ on the canonical unitaries
we see that it has the form $\alpha_{S,\Theta} \beta_x$ for some
$x\in\mathbb{T}^d$, where $\beta$ is the canonical action of
$\mathbb{T}^d$ on $A_\Theta$. Proposition~\ref{P-op} and
Lemma~\ref{L-circle} then yield
$ht(\alpha_{S,-\Theta}) = ht(\alpha^{\op}_{S,-\Theta}) =
ht(\alpha_{S,\Theta})$. By Lemma~\ref{L-tensor} we also have
$ht(\alpha_{S,\Theta}) + ht(\alpha_{S,-\Theta}) \geq
\sum_{| \lambda_i | > 1} \log | \lambda_i |$. The assertion of the theorem
now follows.
\end{proof}

\section{Positive local Voiculescu entropy with respect to products of
canonical unitaries}\label{S-posloc}

Our goal in this section is to obtain some information
concerning positivity of local Voiculescu-Brown entropy with respect
to products of canonical unitaries. We will proceed by first relating
a noncommutative toral automorphism to the corresponding classical
toral automorphism as in the proof of Lemma~\ref{L-tensor} but at a
local level, and then
appealing to a result from \cite{EID} involving local Voiculescu-Brown
entropy in the separable unital commutative setting.
Throughout this section we will denote the canonical unitaries of
the commutative $d$-torus $A_0 \cong C(\mathbb{T}^d )$ by
$f_1 , \dots , f_d$. We will continue to denote the canonical unitaries
of a general noncommutative $d$-torus by $u_1 , \dots , u_d$.

\begin{lemma}\label{L-loctori}
Let $\alpha_{S,\Theta}$ be a noncommutative toral automorphism,
$k_1 , \dots , k_d \in\mathbb{Z}$, and $\lambda\in\mathbb{C}$. Then
$$ ht(\alpha_{S,\Theta} , \{ \lambda
u^{k_1}_1 \cdots u^{k_d}_d \} ) \geq \frac12 ht(\alpha_{S,0} , \{
\lambda f^{k_1}_1 \cdots f^{k_d}_d \} ) .$$
\end{lemma}

\begin{proof}
We may assume that $\lambda = 1$.
Let $C$ be the $\alpha_{S,\Theta} \otimes \alpha_{S,-\Theta}$-invariant
commutative $C^*$-algebra of $A_\Theta \otimes A_{-\Theta}$
identified in the proof of Lemma~\ref{L-tensor}. Denoting by
$v_1 , \cdots , v_d$ the canonical unitaries of $A_{-\Theta}$, we have
\begin{align*}
ht(\alpha_{S,0} , \{ f^{k_1}_1 \cdots f^{k_d}_1 \})
&= ht(\alpha_{S,\Theta} \otimes \alpha_{S,-\Theta} \big| {}_C ,
\{ u^{k_1}_1 \cdots u^{k_d}_d \otimes v^{k_1}_1 \cdots v^{k_d}_d \} ) \\
&\leq ht(\alpha_{S,\Theta} , \{ u^{k_1}_1 \cdots
u^{k_d}_d ) \} ) + ht(\alpha_{S,-\Theta} , \{ v^{k_1}_1 \cdots v^{k_d}_d \} ) ,
\end{align*}
where the last inequality follows from an argument similar to that
in the proof of Proposition 3.10 in \cite{Voi}. As in the proof of
Theorem~\ref{T-half} we identify $A_\Theta$ with
$A^{\op}_{-\Theta}$ via $u_j \mapsto \tilde{v}_j$ and observe that
$\alpha^{\op}_{-\Theta} \in\Aut (A_\Theta )$ has the form
$\alpha_{S,\Theta} \beta_x$ for some $x\in\mathbb{T}^d$, where
$\beta$ is the canonical action of $\mathbb{T}^d$ on $A_\Theta$.
Following Notation~\ref{N-op}, we then have $\widetilde{v^{k_1}_1
\cdots v^{k_d}_d} = \eta \tilde{v}^{k_1}_1 \cdots
\tilde{v}^{k_d}_d$ for some $\eta\in\mathbb{C}$ of unit modulus,
and so
$$ ht (\alpha_{S,-\Theta} , \{ v^{k_1}_1 \cdots v^{k_d}_d \} )
= ht (\alpha^{\op}_{S,-\Theta} , \{ \widetilde{v^{k_1}_1 \cdots
v^{k_d}_d} \} ) = ht(\alpha_{S,\Theta} , \{ u^{k_1}_1 \cdots
u^{k_d}_d \} ) $$ in view of the proofs of Proposition~\ref{P-op}
and Lemma~\ref{L-circle}. The assertion of the lemma now follows
from the above two displays.
\end{proof}

\begin{theorem}\label{T-localpos}
Let $\alpha_{S,\Theta}$ be a noncommutative toral automorphism and
suppose that $S$ has no roots of unity as eigenvalues. Then
$$ ht(\alpha_{S,\Theta} , \{ \lambda u^{k_1}_1
\cdots u^{k_d}_d \} ) > 0 $$
for any non-zero $(k_1 , \dots , k_d )\in\mathbb{Z}^d$ and
non-zero $\lambda\in\mathbb{C}$.
\end{theorem}

\begin{proof}
Since the measure-theoretic toral automorphism associated to
$S$ via Lebesgue measure is ergodic (see \cite{Wal}) and
hence has completely positive
(Kolmogorov-Sinai) entropy \cite{MP}, the topological toral automorphism
associated to $S$ (i.e., the case $\Theta = 0$ at the level of the pure
state space) has completely positive (topological) entropy, i.e., each of
its non-trivial factors has positive topological entropy (see \cite{FPTE}).
Thus $ht(\alpha_{S,0} , \{ x \}) > 0$ for every
non-scalar $x\in A_0 \cong C(\mathbb{T}^d )$ by Corollary
4.4 of \cite{EID}. Lemma~\ref{L-loctori} then yields the result.
\end{proof}

For a general noncommutative toral automorphism
$\alpha_{S,\Theta}$ it follows from Lemma~\ref{L-loctori} that to
conclude that $ht(\alpha_{S,\Theta} , \{ \lambda u^{k_1}_1 \cdots
u^{k_d}_d \} ) > 0$ we need only show that $ht(\alpha_{S,0} , \{
\lambda f^{k_1}_1 \cdots f^{k_d}_d \})
> 0$. If we are simply dealing with a canonical unitary
$u_j$ then this occurs, for example, if the $i$th coordinate axis
in $\mathbb{R}^d$ is not orthogonal to some one-dimensional
subspace of a eigenspace in $\mathbb{R}^d$ corresponding to a real
eigenvalue of $S$ not equal to $\pm 1$. To see this, suppose $L$
is such a one-dimensional subspace and let $\lambda$ be the
associated real eigenvalue. We may assume $|\lambda | > 1$ since
$ht(\alpha^{-1}_{S,\Theta} , \{ u_j \}) = ht(\alpha_{S,\Theta} ,
\{ u_j \})$ (see the proof of Proposition 2.5 in \cite{Br}) and
$\alpha^{-1}_{S,\Theta} = \alpha_{S^{-1},\Theta} \beta_x$ for some
$x\in\mathbb{T}^d$, where $\beta$ is the canonical action of
$\mathbb{T}^d$ on $A_\Theta$. Define the pseudo-metric $d_j$ on
$\mathbb{T}^d$ by
$$ d_j (x,y) = | f_j (x) - f_j (y) | $$
for all $x,y\in\mathbb{T}$, with the unitary $f_i$ being
considered here in the canonical way as a function on the pure
state space $\mathbb{T}^d$. Since the $j$th coordinate subspace of
$\mathbb{R}^d$ is not orthogonal to $L$ there exists a $\delta >
0$ such that, for any $x = (x_1 , \dots , x_d )$ and
$y = (y_1 , \dots , y_d )$ in $L$, the $j$th coordinate distance
$| x_j - y_j |$ is at least
$\delta$ times the Euclidean distance between $x$ and $y$. Now since the
action of $T$ on $L$ is simply multiplication by $\lambda$ it can
be seen via a covering space argument (see the proof of Theorem
24.5 in \cite{DGS}) that there exists a $C>0$ and an $\varepsilon
> 0$ such that for every $n\in\mathbb{N}$ the image of $L$ under
the quotient map onto $\mathbb{T}^d \cong\mathbb{R}^d /
\mathbb{Z}^d$ contains an $(n,\varepsilon )$-separated set of
cardinality at least $C |\lambda |^n$, from which it follows that
the entropy $h_{d_j} (\bar{S})$ is strictly positive, where
$\bar{S}$ is the automorphism of $\mathbb{T}^d$ corresponding to
$S$. Here we are using standard notation and terminology from
topological dynamics (see \cite{DGS}) except that we are allowing
the metric in the definition of entropy to be merely a
pseudo-metric. Now by Theorem 4.3 of \cite{EID} we conclude that
$ht(\alpha_{S,0} , \{ u_j \}) > 0$, as desired.

By a similar argument which allows for the possibility of non-trivial
Jordan cells we have the following more general statement.

\begin{theorem}
Let $\alpha_{S,\Theta}$ be a noncommutative toral automorphism.
Let $j\in \{ 1 , \dots , d \}$ and suppose that the $j$th coordinate
axis in $\mathbb{R}^d$ is not orthogonal to
the span of the generalized eigenspaces associated to the set of real
eigenvalues of $S$ not equal to $\pm 1$. Then
$ht(\alpha_{S,\Theta} , \{ u_j \} ) > 0$.
\end{theorem}

We could evidently furthermore generalize this theorem to handle products of
canonical unitaries, and also formulate a similar result involving pairs
of canonical unitaries and complex eigenvalues not on the unit circle.

\section{Completely positive Voiculescu-Brown entropy}\label{S-K}

In \cite{Nes} Neshveyev showed that, for the von Neumann algebraic
dynamical system arising from a noncommutative toral automorphism
$\alpha_{S,\Theta}$ via
the canonical tracial state on $A_\Theta$, the property of being an
entropic $K$-system (and in particular of having completely positive CNT
entropy), in the case of $S$ having no roots of unity as eigenvalues
(which occurs if and only if $S$ is aperiodic (see \cite{Pet} or
\cite{Wal})), is a consequence of a summability
condition which for $d=2$ is satisfied for a certain countable set of
rotation parameters. We will show in the general unital setting
that completely positive CNT entropy of the von Neumann algebraic dynamical
system arising from a faithful invariant state
implies ``completely positive Voiculescu-Brown entropy''
(i.e., every restriction of the automorphism to a non-trivial
invariant $C^*$-subalgebra has positive Voiculescu-Brown entropy), and that
in the unital case the latter property is equivalent to the positivity of
the unital version of local Voiculescu-Brown entropy with respect to
every non-scalar unitary.

We will use the standard notation for CNT (Connes-Narnhofer-Thirring)
entropy \cite{CNT}.
Let $\alpha$ be an automorphism of a unital $C^*$-algebra $A$ and $\omega$
a faithful $\alpha$-invariant state on $A$. Denoting by $\pi_\omega$ the
GNS representation corresponding to $\omega$, we obtain extensions
$\bar{\alpha}$ and $\bar{\omega}$ of $\alpha$ and $\omega$, respectively,
to $\pi_\omega (A)''$. By definition the automorphism $\bar{\alpha}$ has
completely positive CNT entropy if $h_{\bar{\omega} , \bar{\alpha}} (N) > 0$
for all unital finite-dimensional $C^*$-subalgebras of $\pi_\omega (A)''$
which are different from the scalars. Entropic $K$-systems have completely
positive CNT entropy, and the two notions coincide in the
commutative case. For definitions and discussions see \cite{QKS,GN}. Note
that in \cite{GN} what we are referring to for clarity as
``completely positive CNT entropy'' is simply called
``completely positive entropy.''

\begin{definition}
An automorphism $\alpha$ of a non-trivial exact $C^*$-algebra $A$ is
said to have {\it completely positive Voiculescu-Brown
entropy} if $ht(\alpha\big| {}_B ) > 0$ for every non-zero
$\alpha$-invariant $C^*$-subalgebra $B\subseteq A$ which, if $A$
is unital, is not equal to the scalars.
\end{definition}

\begin{remark}
If $A$ is unital then in the above definition we may take the
$C^*$-subalgebras $B$ to be unital, for if $B\subseteq A$ is an
$\alpha$-invariant $C^*$-subalgebra not containing the unit of $A$ then
the Voiculescu-Brown entropies of the restrictions $\alpha\big|
{}_B$ and $\alpha\big| {}_{B+\mathbb{C}1}$ agree
by Lemma 1.7 of \cite{Br}. Thus in the
separable unital commutative situation we recover the
topological-dynamical notion of completely positive entropy, which refers
to the absence of non-trivial factors with zero topological
entropy \cite{FPTE}.
\end{remark}

\begin{proposition}\label{P-cpe}
Let $\alpha$ be an automorphism of a unital exact $C^*$-algebra $A$
and $\omega$ a faithful $\alpha$-invariant state on $A$, and suppose
that the extension $\bar{\alpha}$ of $\alpha$ to $\pi_\omega (A)''$
has completely positive CNT entropy. Then $\alpha$ has completely positive
Voiculescu-Brown entropy.
\end{proposition}

\begin{proof}
Let $B\subseteq A$ be a unital $\alpha$-invariant $C^*$-subalgebra
different from the scalars. By choosing a non-scalar element in
$B$ and taking its real part (or imaginary part if its real part
is a scalar) we obtain a self-adjoint element $b\in B$ whose
spectrum contains more than one point, and we can construct an
injective unital $^*$-homomorphism $\gamma : \mathbb{C}^2 \to
\pi_\tau (C(\spec (b)))'' \subseteq \pi_\tau (B)''$. Since
$\bar{\alpha}$ has completely positive CNT entropy we have
$h_{\bar{\omega} , \bar{\alpha}} (\gamma ) > 0$. For any $\delta >
0$ it is easy to construct, using a suitable two-element partition
of unity in $C(\spec (b))$, a unital positive map $\gamma' :
\mathbb{C}^2 \to C(\spec (b))\subseteq B$ (which is necessarily
completely positive since $\mathbb{C}^2$ is commutative) such that
$\| \gamma - \gamma' \|_\tau < \delta$, where
$$ \| \gamma - \gamma' \|_\tau =
\sup_{\{ x\in\mathbb{C}^2 \, : \, \| x \| \leq 1 \}} \| \gamma (x)
- \gamma' (x) \|_\tau $$ and $\| a \|_\tau = ( \tau (a^*
a))^{1/2}$ for all $a\in A$. It follows from Theorem VI.3 of
\cite{CNT} that if $\delta$ is chosen small enough then
$$ h_{\omega|_B , \alpha|_B } (\gamma' ) \geq
h_{\bar{\omega} , \bar{\alpha}} (\gamma' ) > 0 $$
with the first inequality being immediate from the definition of CNT
entropy. Proposition 9 of \cite{Dyk} then yields
$$ ht(\alpha\big| {}_B ) \geq h_{\omega|_B} (\alpha\big| {}_B )
> 0 , $$
as desired.
\end{proof}

Recall from the beginning of Section~\ref{S-K} that
$ht' (\alpha , \Omega )$ denotes the unital version of the local
Voiculescu-Brown entropy $ht(\alpha , \Omega )$.
We will next show that, in the unital case,
completely positive Voiculescu-Brown entropy is equivalent to positivity
of the unital version of local
Voiculescu-Brown entropy with respect to every non-scalar unitary.
For this we will need
the following Kolmogorov-Sinai-type property, which is similar to that of
Lemma A.2 in \cite{Oz}.

\begin{lemma}\label{L-KS}
Let $\alpha$ be an automorphism of a unital exact $C^*$-algebra
$A$. If $\Omega_1 \subseteq \Omega_2 \subseteq \Omega_3
\subseteq\cdots$ is a nested sequence of finite sets of unitaries
in $A$ such that
$\bigcup_{k\in\mathbb{N},n\in\mathbb{Z}} \alpha^n
\Omega_k$ generates $A$ as a $C^*$-algebra, then
$$ ht(\alpha ) = \sup_{k\in\mathbb{N}} ht' (\alpha , \Omega_k ) . $$
\end{lemma}

\begin{proof}
Given a unital completely positive map $\varphi$ from $A$ into any unital
$C^*$-algebra $B$, Lemma 3.1 of \cite{JOR} yields
$$ \| \varphi (x^* y) - \varphi (x)^* \varphi (y) \| \leq \| \varphi (x^* x)
- \varphi (x)^* \varphi (x) \|^{1/2}\| \varphi (y^* y) -
\varphi (y)^* \varphi (y) \|^{1/2} $$
for all $x,y\in A$, and so in particular for any unitaries $u,v \in A$
we have
\begin{align*}
\| \varphi (uv) - \varphi (u) \varphi (v) \| &\leq \| 1 - \varphi (u)^*
\varphi (u) \|^{1/2} \| 1 - \varphi (v)^* \varphi (v) \|^{1/2} \\
&\leq (\| u^* u - \varphi (u)^* u \| + \| \varphi (u)^* u -
\varphi (u)^* \varphi (u) \| )^{1/2} \\
&\hspace*{8mm}\ \times (\| v^* v - \varphi (v)^* v \| +
\| \varphi (v)^* v - \varphi (v)^* \varphi (v) \|)^{1/2} \\
&\leq 2 \| \varphi (u) - u \|^{1/2} \| \varphi (v) - v \|^{1/2}
\end{align*}
whence
\begin{align*}
\| \varphi (uv) - uv \| &\leq \| \varphi (uv) - \varphi (u) \varphi (v) \| +
\| \varphi (u) \varphi (v) - u\varphi (v) \| + \| u\varphi (v) - uv \| \\
&\leq (\| \varphi (u) - u \|^{1/2} + \| \varphi (v) - v \|^{1/2} )^2 .
\end{align*}
We can now proceed along the lines of
the proofs of Propositions 1.4 and 3.4 of \cite{Voi} to obtain the result.
\end{proof}

\begin{proposition}\label{P-cpeunit}
Let $\alpha$ be an automorphism of a unital exact $C^*$-algebra
$A$. Then $\alpha$ has completely positive Voiculescu-Brown
entropy if and only if $ht' (\alpha , \{ u \} ) > 0$ for every
non-scalar unitary $u\in A$.
\end{proposition}

\begin{proof}
For the ``only if'' direction we can consider for any non-scalar
unitary $u\in A$ the unital $\alpha$-invariant $C^*$-subalgebra it generates
and appeal to Lemma~\ref{L-KS}. For the ``if'' direction, we need
simply observe that
every unital $C^*$-algebra different from the scalars contains a non-scalar
unitary, as can be obtained by applying the functional calculus to the
real part (or imaginary part if the real part is a scalar) of any
non-scalar element.
\end{proof}

For a real skew-symmetric $d\times d$ matrix $\Theta$ the
noncommutative $d$-torus $A_\Theta$ may alternatively be described
as the universal unital $C^*$-algebra generated by unitaries $\{
u_g \}_{g\in \mathbb{Z}^d}$ subject to the relations
$$ u_g u_h = \beta (g,h) u_{g+h}, $$
where $\beta : \mathbb{Z}^d \times\mathbb{Z}^d \to\mathbb{T}$ is a bicharacter
satisfying
$$ \beta (g,h) \beta (h,g)^{-1} = e^{2\pi ig\cdot\Theta h}. $$

\begin{theorem}\label{T-sum}
Let $\alpha_{S,\Theta}$ be a noncommutative toral automorphism,
and suppose that $S$ is $\beta$-preserving and has no roots of unity as
eigenvalues, and that
$$ \sum_{n\in\mathbb{Z}} | 1 - \beta (g,S^n h ) | < \infty $$
for all $g,h\in \mathbb{Z}^d$. Then $\alpha_{S,\Theta}$ has
completely positive Voiculescu-Brown entropy and
$ht' (\alpha_{S,\Theta} , \{ u \} ) > 0$ for every non-scalar
unitary $u\in A$.
\end{theorem}

\begin{proof}
By Theorem 2 of \cite{Nes} the hypotheses of the theorem statement
imply that the von Neumann algebraic system obtained from
$\alpha_{S,\Theta}$ via the canonical tracial state on $A_\Theta$
is an entropic $K$-system, and hence has completely positive CNT
entropy. Propositions~\ref{P-cpe} and \ref{P-cpeunit} then yield
the desired conclusion.
\end{proof}

Theorem~\ref{T-sum} applies in particular in the case $d=2$ when the matrix
$S$ has eigenvalues off the unit circle and
the rotation parameter $\theta$ of the rotation $C^*$-algebra $A_\theta$
($=A_\Theta$ for $\Theta = (\theta_{jk})_{j,k=1,2}$ with
$\theta_{12} = \theta$) lies in $\mathbb{Z} + 2\mathbb{Z} \lambda^2$
where $\lambda$ is the (necessarily real) eigenvalue of $S$ of largest
absolute value (see \cite{Nes}).

Finally, we would like to point out that the argument in \cite{Nes}
is not quite complete. Indeed in the proof of the lemma in \cite{Nes}
it is incorrectly taken to be the case that a matrix
$T\in GL(n,\mathbb{Z})$ is aperiodic if and only if it has no eigenvalues
on the unit circle. (Aperiodicity is defined as the non-existence of
non-trivial finite orbits of $T$ acting on $\mathbb{Z}^n$ and is equivalent
to $T$ having no roots of unity as eigenvalues and also equivalent to the
ergodicity of the measure-theoretic automorphism of $\mathbb{T}^n$
associated to $T$ via Lebesgue measure (see \cite{Pet} or \cite{Wal}, and
also Section 24 of \cite{DGS}).)
However, if we let $\mathbb{R}^n = V_1 \oplus V_2$
be the decomposition of $\mathbb{R}^n$ corresponding to the
eigenvalues of $T$ of modulus at least one and strictly less than one,
respectively, and denote by $P_1$ and $P_2$ the associated projections,
then the proof of the lemma in \cite{Nes} demonstrates that, for large
$n$, if $y_1 + T^n y_2 + \dots + T^{n(k-1)} y_k = 0$ then
$P_2 (y_1 ) = \dots = P_2 (y_k ) = 0$. The argument in the first
paragraph of p.\ 191 of \cite{Katz} then shows that $P_2$ is
injective on $\mathbb{Z}^n$, whence $y_1 = \dots = y_k = 0$, as desired.

\end{document}